\documentclass[12pt]{article}
\usepackage{amsmath,amsthm,amsfonts}

\newtheorem{thm}{Theorem}[section]
\newtheorem{lemma}[thm]{Lemma}

\newtheorem{cor}[thm]{Corollary}
\newtheorem{defin}[thm]{Definition}
\newtheorem{rem}[thm]{Remark}

\newtheorem{problem}[thm]{Problem}
\newtheorem{question}[thm]{Question}

\newcommand{\R}{{\mathbb{R}}}

\newcommand{\T}{{\mathbb{T}}}
\newcommand{\Z}{{\mathbb{Z}}}
\newcommand{\N}{{\mathbb{N}}}
\newcommand{\C}{{\mathbb{C}}}

\newcommand{\cA}{{\mathcal{A}}}
\newcommand{\cD}{{\mathcal{D}}}

\newcommand{\cF}{{\mathcal{F}}}
\newcommand{\cG}{{\mathcal{G}}}

\newcommand{\cS}{{\mathcal{S}}}

\newcommand{\cC}{{\mathcal{C
}}}

\newcommand{\Cal}{{\hbox{\it Cal\,}}}

\newcommand{\Ham}{{\hbox{\it Ham\,}}}

\newcommand{\id}{{\text{{\bf 1}}}}
\newcommand{\tHam}{\widetilde{\hbox{\it Ham}\, }}
\newcommand{\Symp}{{\hbox{\it Symp} }}
\newcommand{\supp}{{supp}}

\newcommand{\Qed}{\hfill \qedsymbol \medskip}

\begin{document}

\title{Quasi-states and symplectic intersections\\
}

\renewcommand{\thefootnote}{\alph{footnote}}

\author{\textsc Michael Entov$^{a}$\ and Leonid
Polterovich$^{b}$}

\footnotetext[1]{Partially supported by the Technion President
Fund, by the Israel Science Foundation grant $\#$ 68/02 and by M.
and C.Papo Research Fund.} \footnotetext[2]{Partially supported by
the Israel Science Foundation grant $\#$ 11/03.}

\date{\today}

\maketitle

\begin{abstract}
\noindent We establish a link between symplectic topology and a
recently emerg\-ed branch of functional analysis called the theory
of quasi-states and quasi-measures (also known as topological
measures). In the symplectic context quasi-states can be viewed as
an algebraic way of packaging certain information contained in
Floer theory, and in particular in spectral invariants of
Hamiltonian diffeomorphisms introduced recently by Yong-Geun Oh.
As a consequence we prove a number of new results on rigidity of
intersections in symplectic manifolds. This work is a part of a
joint project with Paul Biran.
\end{abstract}

\vfil \eject

\tableofcontents

\renewcommand{\thefootnote}{\arabic{footnote}}
\vfil \eject

\section{Introduction}
\label{sec-intro}

Rigidity of intersections is a class of phenomena in symplectic
topology meaning that certain subsets of a symplectic manifold
intersect each other in more points than dictated by algebraic and
differential topology (see \cite{Biran} for an excellent survey).
In this paper we show that such rigidity phenomena in a closed
symplectic manifold $M$ sometimes formally follow from the
existence of real-valued functionals with some interesting
algebraic properties on the Poisson algebra $C^\infty (M)$.

On the one hand, these functionals are related to the notions of
{\it quasi-state} and {\it quasi-measure} (which have been
recently called {\it topological measures}) on $M$ (see
Section~\ref{sec-quasi-objects}) which originate in quantum
mechanics \cite{Aar69},\cite{Aar70} and have been a subject of
intensive study in recent years following the paper \cite{Aar91}
by J.F.Aarnes.

On the other hand, they are linked to a group-theoretic notion of
{\it quasi-morphism} (see e.g. \cite{Kotschick}) which already
appeared in the context of symplectic topology in \cite{EP},
\cite{BEP}. The symplectic quasi-states on the Poisson-Lie algebra
of functions on certain symplectic manifolds $M$ considered below
arise as an infinitesimal version of the Calabi quasi-morphism
introduced in \cite{EP}. This quasi-morphism is defined on the
universal cover $\tHam (M)$ of the group $\Ham (M)$ of Hamiltonian
diffeomorphisms of $M$.

All the above-mentioned functionals are constructed by means of
Floer theory for Hamiltonian flows on $M$ and can be viewed as an
algebraic way of packaging certain information contained in that
theory.

Throughout the paper $M$ always stands for a closed connected
symplectic manifold with a symplectic form $\omega$. For technical
reasons we assume that $M$ is {\bf rational}, i.e. the image of
$\pi_2 (M)$ under the cohomology class of $\omega$ is a discrete
subgroup of $\R$. Furthermore, we assume that $M$ is {\bf strongly
semi-positive}, that is
\begin{equation} \label{eqn-strongly-semi-pos} 2-n\leq c_1 (A)<
0 \Longrightarrow \omega (A)\leq 0, \ {\rm for\ any}\ A\in \pi_2
(M),
\end{equation}
where $c_1$ stands for the 1st Chern class of $(M,\omega)$. For
instance, every symplectic $4$-manifold is strongly semi-positive.
Another interesting class of examples is given by {\it spherically
monotone} symplectic manifolds, which means that $\left. [\omega]
\right|_{\pi_2 (M)}$ is a positive multiple of $\left.
c_1\right|_{\pi_2 (M)}$. Note that this condition automatically
implies strong semi-positivity and rationality of $M$.

\medskip
\noindent {\sc Organization of the paper:} The next section
contains our main results on symplectic intersections. In Section
3 we focus on a special class of symplectic manifolds $M$ which,
for instance, includes monotone products of complex projective
spaces. After a brief review of quasi-states and quasi-measures,
we introduce symplectic quasi-states on the algebra $C(M)$ which
turn out to be useful for symplectic intersections in $M$. In
Section 4 we present a weaker notion of a partial symplectic
quasi-state and its applications to non-displaceability phenomenon
on more general symplectic manifolds. In Section 5 we review
spectral invariants of Hamiltonian diffeomorphisms introduced
recently by Y.-G.~Oh. In Sections 6 and 7 these invariants are
used in order to construct the above-mentioned (partial)
symplectic quasi-states. In Section 8 we discuss symplectic
quasi-states on surfaces. The reader will see that some innocently
looking basic questions in this direction require more advanced
tools of the theory of quasi-states and quasi-measures. Section 9
contains some applications (in the spirit of our paper \cite{BEP}
with P.~Biran) of our results to the Lagrangian intersection
problem. In Section 10 we discuss the history and the physical
meaning of quasi-states. In addition, in Sections 8-10 we present
a number of open problems.

\section{Results on symplectic intersections}
\label{sec-sample}

We say that a subset $X$ of $M$ is {\it displaceable} if there
exists a Hamiltonian diffeomorphism $\phi \in \Ham(M)$ so that
$$\phi(X) \cap \text{Closure}(X) = \emptyset.$$ Otherwise, we call
$X$ {\it non-displaceable}. For instance, an open hemisphere in
$S^2$ with the Euclidean area form is displaceable while the
closed hemisphere is not.

A linear subspace $\cA \subset C^{\infty}(M)$ is called {\it
Poisson-commutative}, if $\{F,G\}=0$ for all $F,G \in \cA$, where
$\{\cdot,\cdot\}$ stand for the Poisson brackets. Given a
finite-di\-men\-sio\-nal Poisson-commutative subspace $\cA\subset
C^\infty (M)$, its {\it moment map} $\Phi_\cA: M\to \cA^*$ is
defined as
\[
\langle \Phi_\cA (x), F\rangle = F (x).
\]
Non-empty subsets of the form $\Phi_\cA^{-1} (p), \; p \in \cA^*$,
are called {\it fibers} of $\cA$.

\begin{thm}
\label{thm-commut-alg-has-nondispl-fiber} Any
fi\-nite-di\-men\-sio\-nal Poisson-com\-mu\-ta\-tive sub\-space of
\break $C^\infty (M)$ has at least one non-displaceable fiber.
Moreover, if every fiber has a finite number of connected
components, there exists a fiber with a non-displaceable connected
component.
\end{thm}

Poisson-commutative subspaces naturally appear when $M$ is
equipped with the structure of a (singular) Lagrangian fibration.
In this case Theorem~\ref{thm-commut-alg-has-nondispl-fiber} shows
that the fibration has at least one non-displaceable fiber. For
instance, we have the following corollary, where the fibration is
given by the moment map of a Hamiltonian torus action.

\noindent
\begin{cor}
\label{cor-toric} Assume that $M^{2n}$ is equipped with a
Hamiltonian action of $\T^n$. Then at least one Lagrangian orbit
of this action is non-displaceable.
\end{cor}

\medskip
\noindent {\bf Proof of the Corollary:} Let $\cA$ be the span of
the coordinate functions associated to the moment map of the
action. Every fiber of $\cA$ is a fiber of the moment map: it is
either a Lagrangian torus, or an isotropic torus of dimension less
than $n$. The latter are displaceable (see e.g.
\cite{Bir-Ciel-CMH}). Hence the result follows immediately from
Theorem~\ref{thm-commut-alg-has-nondispl-fiber}. \Qed

\noindent
\begin{defin} \label{def-stem}
{\rm A closed subset $X\subset M$ is called a {\it stem}, if there
exists a finite-dimensional Poisson-commutative subspace
$\cA\subset C^\infty (M)$ so that $X$ is a fiber of $\cA$ and each
fiber of $\cA$, other than $X$, is displaceable.}
\end{defin}

\medskip
\noindent Note that the image of a stem under {\it any}
symplectomorphism of $M$ is again a stem.

\medskip
\noindent Theorem~\ref{thm-commut-alg-has-nondispl-fiber}
guarantees that {\it every stem is non-displaceable}. This result
can be strengthened for a special class of symplectic manifolds as
follows.

\begin{thm}
\label{thm-genuine-quasi-objects} Suppose $M$ is one of the
following symplectic manifolds: \noindent $\C P^n$, a complex
Grassmannian, $\C P^{n_1}\times\ldots\times\C P^{n_k}$ with a
monotone product symplectic structure, the monotone symplectic
blow-up of $\C P^2$ at one point. Then any two stems in $M$ have
an non-empty intersection.
\end{thm}

In particular, {\it a stem in such an $M$ cannot be displaced from
itself by any (not necessarily Hamiltonian) symplectomorphism}.

Here is a sample corollary of this theorem. Consider the 2-sphere
$S^2$ with a symplectic form $\omega$ of total area 1. Define a
class $\cG_{S^2}$ of closed subsets $\Gamma \subset S^2$ with the
following property: The complement $S^2 \setminus \Gamma$ has a
finite number of connected components, and each of them is
homeomorphic to a disc and has area $\leq \frac{1}{2}$. For
instance, one can take an equator, or the 1-skeleton of a
piecewise smooth triangulation of $S^2$ with small enough
2-dimensional faces.

\noindent
\begin{cor}
\label{thm-skeleta} Let $M$ be the direct product of $m$ copies of
$(S^2, \omega)$ and let $\Gamma_i, \Gamma_i^\prime\in \cG_{S^2}$,
$i=1,\ldots, m$. Then the subsets $\Gamma_1\times\ldots\times
\Gamma_m$ and $\phi(\Gamma_1^\prime\times\ldots \times
\Gamma_m^\prime)$ have a non-empty intersection for every
symplectomorphism $\phi$ of $M$.
\end{cor}

\medskip
\noindent{\bf Proof of the Corollary:} Note that a direct product
of stems is a stem. Hence it suffices to verify that every $\Gamma
\in \cG_{S^2}$ is a stem. Let $U_1,...,U_d$ be the connected
components of $S^2 \setminus \Gamma$.

Take smooth functions $H_1,...,H_d$ as follows: $H_i$ vanishes on
$S^2 \setminus U_i$ and $H_i$ is strictly positive on $U_i$. The
existence of such $H_1,...,H_d$ follows easily from the fact that
any closed subset of $\R^2$ is the zero-level set of some smooth
real-valued function on $\R^2$ (see e.g. \cite{Narasimhan}, Lemma
1.4.13).

Put $\cA = \text{Span}_{\R}(H_1,...,H_d)$. Clearly $\cA$ is
Poisson-commutative and $\Gamma = \Phi_\cA^{-1}(0)$ is its fiber.
All other fibers are closed subsets of one of the $U_i$'s, and
hence are displaceable. Therefore $\Gamma$ is a stem and the
result follows from Theorem~\ref{thm-genuine-quasi-objects}. \Qed

\bigskip
Here is another corollary of
Theorem~\ref{thm-commut-alg-has-nondispl-fiber}. Let $\T^2$ be a
torus with coordinates $p,q\in \R/\Z$ and the symplectic form $dp
\wedge d q$. Equip $M\times \T^2$ with the product symplectic
structure and assume that the resulting symplectic manifold is
strongly semi-positive and rational. Denote by $S$ a meridian
$p=const$ of $\T^2$.

\begin{cor}
Assume $X\subset M$ is a stem. Then $X\times S\subset M\times
\T^2$ is non-displaceable.
\end{cor}

\medskip
\noindent{\bf Proof of the Corollary:} Let $\cA\subset C^\infty
(M)$ be a finite-dimensional Poisson-commutative subspace such
that the stem $X$ is its only non-displaceable fiber. Lift to
$M\times \T^2$ the functions on $M$ that belong to $\cA$ as well
as the functions $\sin 2\pi p, \cos 2\pi p$ on $\T^2$. All these
lifts together span a Poisson-commutative subspace
$\cA^\prime\subset C^\infty (M\times \T^2)$ such that each of its
fibers is a direct product of a fiber of $\cA$ and a meridian of
$\T^2$.

Theorem~\ref{thm-commut-alg-has-nondispl-fiber} says that
$\cA^\prime$ must have a non-displaceable fiber $Y$. Since $X$ is
the only non-displaceable fiber of $\cA$, the fiber $Y$ has to
have the form $Y= X\times S^\prime$ for some meridian $S^\prime$
of $\T^2$. But any two meridians of $\T^2$ can be mapped into each
other by a symplectomorphism of $\T^2$ -- hence the products of
these meridians with $X$ can be mapped into each other by a
symplectomorphism of $M\times \T^2$. Thus if $X\times S^\prime$ is
non-displaceable, then $X\times S$ has to be non-displaceable as
well. \Qed

\section{Quasi-states and quasi-measures}
\label{sec-quasi-objects}

Write $C(M)$ for the commutative (with respect to multiplication)
Banach algebra of all continuous functions on $M$ endowed with the
uniform norm. For a function $F\in C(M)$ denote by $\cA_F$ the
uniform closure of the set of functions of the form $p \circ F$,
where $p$ is a real polynomial. A (not necessarily linear)
functional $\zeta: C(M) \to \R$ is called {\it a quasi-state}
\cite{Aar91}, if it satisfies the following axioms:

\medskip
\noindent {\bf Quasi-linearity:} $\zeta$ is linear on $ \cA_F$ for
every $F \in C(M)$ (in particular $\zeta$ is homogeneous);

\medskip
\noindent {\bf Monotonicity:} $\zeta (F) \leq \zeta (G)$ for $F
\leq G$;

\medskip
\noindent {\bf Normalization:} $\zeta (1) = 1.$

A quasi-state is called {\it symplectic}, if it has the following
additional properties:

\medskip
\noindent {\bf Strong quasi-additivity:} $\zeta (F+G)=\zeta
(F)+\zeta (G)$ for all smooth functions $F,G$ which commute with
respect to the Poisson bracket: $\{F,G\} = 0$;

\medskip
\noindent{\bf Vanishing:} $\zeta (F)= 0$, provided $\supp\, F$ is
displaceable;

\medskip
\noindent {\bf Symplectic invariance:} $\zeta (F) = \zeta (F\circ
f)$ for every symplectic diffeomorphism $f \in \Symp_0\, (M)$
(here $\Symp_0 (M)$ stands for the identity component of the group
$\Symp\, (M)$ of symplectomorphisms).

\medskip
\noindent Note that strong quasi-additivity together with
homogeneity yields quasi-linearity. Indeed, if $F$ is smooth, then
$\{F,p\circ F\}=0$ for every polynomial $p$. Observing that
$\zeta$ is continuous in the uniform topology because of the
monotonicity and normalization axioms, one can easily extend the
result for a general continuous $F$.

\medskip
\noindent
\begin{thm} \label{thm-sympl-quasistate}
Suppose $M$ is one of the following symplectic manifolds:
\noindent $\C P^n$, a complex Grassmannian, $\C
P^{n_1}\times\ldots\times\C P^{n_k}$ with a monotone product
symplectic structure, the monotone symplectic blow-up of $\C P^2$
at one point. Then $C(M)$ admits a symplectic quasi-state.
\end{thm}

In \cite{Aar91} Aarnes proved a generalized Riesz representation
theorem which associates to each quasi-state $\zeta$ a {\it
quasi-measure} $\tau_\zeta$, that is a "measure" which is finitely
additive but not necessarily sub-additive. More precisely, denote
by $\cS$ the collection of all subsets of $M$ which are {\bf
either open or closed}. A {\it quasi-measure} (recently called a
{\it topological measure} in the literature) on $M$ is a
$[0,1]$-valued set-function $\tau$ on $\cS$ such that

\smallskip
1) $\tau (M) = 1$;

\smallskip
2) $X_1\subset X_2$ $\Rightarrow$ $\tau (X_1) \leq \tau (X_2)$ for
all $X_1,X_2 \in \cS$;

\smallskip
3) $\tau (X_1\sqcup\ldots \sqcup X_k) = \tau
(X_1)+\ldots+\tau(X_k)$ for all $X_1,...,X_k \in \cS$ with
$X_1\sqcup\ldots \sqcup X_k \in \cS$;

\smallskip
4) For every open subset $X$ one has $\tau (X) = \sup\tau(A)$,
where the supremum is taken over all closed subsets $A \subset X$.

\smallskip
\noindent The relation between a quasi-state $\zeta$ and the
corresponding quasi-measure $\tau_\zeta$ is the following
\cite{Aar91}. Given a closed $X\subset M$, consider the set
$\cF_X$ of smooth functions $M\to [0,1]$ which are identically
equal to $1$ on $X$. A quasi-state $\zeta$ is bounded on $\cF_X$
by $0$ and $1$ and therefore one can define
\begin{equation}
\label{eqn-tau-zeta} \tau_\zeta (X) := \inf_{F\in \cF_X} \zeta
(F).
\end{equation}
Intuitively, $\tau_\zeta (X)$ is the "value" of the functional
$\zeta$ on the (discontinuous) characteristic function of $X$. For
an open subset $Y$ put $\tau_\zeta (Y) = 1-\tau_\zeta (M \setminus
Y).$

\begin{lemma} \label{lem-quasim-1} Assume a closed connected symplectic
manifold $M$ admits a symplectic quasi-state $\zeta$. Denote by
$\tau$ the corresponding quasi-measure. Then $\tau (X) = 1$ for
every stem $X \subset M$.
\end{lemma}

\medskip
\noindent {\bf Proof:} Let $\cA \subset C^{\infty}(M)$ be a
finitely generated Poisson-commutative subspace. Denote by $\Delta
\subset \cA^*$ the image of the moment map $\Phi_{\cA}$. Write
$C^{\infty}_0(\cA^*)$ for the space of all smooth compactly
supported functions on $\cA^*$. Note that the functional $$I:
C^{\infty}_0(\cA^*) \to \R,\; G \mapsto \zeta(\Phi_{\cA}^* G),$$
is a positive distribution\footnote{Recall that a distribution
(that is a continuous linear functional) on $C^{\infty}_0 (\R^N)$
is called {\it positive } if it takes non-negative values on
non-negative functions.} (use the strong quasi-additivity and the
monotonicity axioms of $\zeta$). Hence it defines a measure
$\sigma$ on $\cA^*$ so that $I (G) = \int_{\cA^*} G d\sigma$ (see
e.g. \cite{Gelfand-Vilenkin}, Ch. 2, Sec. 2). By the normalization
axiom, $\sigma$ is a probability measure. Obviously, $\supp\,
\sigma\subset \Delta$. The vanishing axiom yields that if
$\Phi_{\cA}^{-1}(p)$ is displaceable for some $p \in \Delta$, then
$p \notin \supp\, \sigma$. Thus, if $X = \Phi_{\cA}^{-1}(p_0)$ is
a stem associated to $\cA$, the measure $\sigma$ must be the Dirac
measure at $p_0$. Using this and considering in the definition of
$\tau (X)$ the functions $F\in\cF_X$ of the form $F = \Phi_{\cA}^*
G$, $G\in C^{\infty}_0(\cA^*)$, one readily gets $\tau (X) = 1$.
\Qed

The proof of the lemma shows that if $\tau$ is a quasi-measure
defined by a symplectic quasi-state, and $\cA \subset
C^{\infty}(M)$ is a finitely generated Poisson-commutative
subspace, the push-forward of $\tau$ by the moment map
$\Phi_{\cA}$ is a {\it genuine measure} on the image of
$\Phi_{\cA}$. In case when $\tau$ comes from a quasi-state which
is not strongly quasi-additive (and thus not symplectic), this may
no longer be true and moreover such a quasi-measure may vanish on
a stem -- see Remark~\ref{rem-non-sympl-qstates-higher-dim}.

\medskip
\noindent {\bf Proof of Theorem~\ref{thm-genuine-quasi-objects}
(assuming Theorem~\ref{thm-sympl-quasistate}):} According to
Theorem~\ref{thm-sympl-quasistate}, any $M$ mentioned in the
hypothesis of Theorem~\ref{thm-genuine-quasi-objects} admits a
symplectic quasi-state. Let $\tau$ be the corresponding
quasi-measure and let $X,Y\subset M$ be stems.
Lemma~\ref{lem-quasim-1} implies that $\tau (X) = \tau (Y) = 1$.
If $X$ and $Y$ do not intersect, we have $\tau(X \cup Y) = \tau(X)
+ \tau (Y) = 1+ 1 =2,$ and we get a contradiction with $\tau(X
\cup Y) \leq \tau(M) = 1$. \Qed

\section{What happens on more general symplectic manifolds?}

Let $\zeta: C(M) \to \R$ be a (not necessarily quasi-linear)
functional which satisfies monotonicity, normalization, vanishing
and invariance axioms from the previous section. Assume that it
has two additional properties:

\medskip
\noindent {\bf Partial additivity:} If $F_1, F_2\in C^\infty (M)$,
$\{ F_1, F_2\} =0$ and the support of $F_2$ is displaceable, then
$\zeta (F_1 + F_2) = \zeta (F_1)$;

\medskip
\noindent {\bf Semi-homogeneity:} $\zeta (\lambda F) = \lambda
\zeta (F)$ for any $F$ and any $\lambda \in \R_{\geq 0}$.

\medskip
\noindent We call $\zeta$ {\it a partial symplectic quasi-state}.

\begin{thm} \label{thm-part-quas}
Let $M$ be a strongly semi-positive and rational closed connected
symplectic manifold. Then $C(M)$ admits a partial symplectic
quasi-state.
\end{thm}

\noindent Theorem~\ref{thm-part-quas} will be proved in
Section~\ref{sec-existence-part-qstate}.

\medskip
\noindent {\bf Proof of
Theorem~\ref{thm-commut-alg-has-nondispl-fiber} (assuming
Theorem~\ref{thm-part-quas}):} Let $\zeta$ be a partial symplectic
quasi-state. Assume on the contrary that all fibers of $\cA$ are
displaceable. Choose an open covering $U:=\{U_1,...,U_d\}$ of the
image $\Delta$ of the moment map $\Phi_{\cA}$ so that the
preimages $\Phi^{-1}(U_i)$ are displaceable. Let
$\rho_1,...,\rho_d$ be a partition of unity associated to $U$,
that is $\supp\, \rho_i \subset U_i$ and $\sum_{i=1}^d \rho_i
\Big{|}_\Delta = 1$. Note that $\zeta(\Phi^*\rho_i)=0$ by
vanishing property. Using the normalization and the partial
additivity, we get $$1= \zeta(1) = \zeta \Big{(}\sum_{i=1}^d
\Phi^*\rho_i\Big{)} = \sum_{i=1}^d \zeta(\Phi^*\rho_i)=0\;,$$ and
we get a contradiction.

A similar argument shows that, if any fiber of $\cA$ has a finite
number of connected components, then at least one connected
component of some fiber of $\cA$ has to be non-displaceable.\Qed

\section{Spectral numbers -- review}
\label{sec-spectral-numbers}

We review a few basic facts about the spectral numbers of
Hamiltonian dif\-feo\-morphisms introduced by Yong-Geun Oh
\cite{Oh-spectral} (see also \cite{Viterbo, Schwarz} for earlier
versions of this theory). For the precise definitions and further
details see \cite{Oh-spectral}, \cite{EP} and
\cite{McD-Sal-pshc-book}. We assume here that $M$ is strongly
semi-positive and rational. The strong semi-positivity of $M$ is
needed to guarantee that the moduli spaces of pseudo-holomorphic
curves involved in the definitions of Floer and quantum homology
and the isomorphism between them are well-behaved. In view of the
developments \cite{Fu-Ono}, \cite{Liu-Ti}, \cite{Liu-Ti1},
\cite{Lu}, concerning Floer theory for general symplectic
manifolds, it is likely that the strong semi-positivity of $M$ is
not essential for the existence of spectral numbers. The
assumption that $M$ is rational is needed to guarantee the
spectrality property below\footnote{For the same reason the
rationality assumption should be added to the results 2.5.3,
2.5.4, 2.6.1 and to part 4 of 2.4.2 in \cite{En-comm} which
involve the spectral numbers.}, though it is likely that
eventually this assumption will also be removed, see
\cite{Oh-last}.

By ${\it \overline{spec}}\, (H)$ we denote the action spectrum of
a Hamiltonian $H$. Recall that it is the set of critical values of
the action functional defined by $H$ on the universal cover of the
space of free contractible loops in $M$.

A time-dependent Hamiltonian $H: M\times S^1\to\R$ is called {\it
normalized} if
$$
\int_M H(\cdot,t) \omega^n =0 \;\;\text{for}\;\text{all}\; t \in
S^1.
$$
It turns out that ${\it \overline{spec}}\, (H_1) = {\it
\overline{spec}}\, (H_2)$ for any normalized $H_1, H_2$ generating
the same element $\phi\in\tHam (M)$. Thus one can define ${\it
spec}\, (\phi)$ for any $\phi\in\tHam (M)$ as ${\it
\overline{spec}}\, (H)$ for any normalized $H$ generating $\phi$.

Denote by $QH_* (M)$ the quantum homology ring of $M$ (with
coefficients in $\C$) and by $*$ the product in that ring. The
fundamental class $[M]$ is the unit in the ring. To each non-zero
quantum homology class $a\in QH_* (M)$ and each time-dependent
Hamiltonian $H: M\times S^1\to\R$ one can associate a {\it
spectral number} $\bar{c} (a, H)$. Spectral numbers have the
following properties which are relevant for us:

\begin{description}

\item[{\bf (Spectrality)}]\ $\bar{c} (a, H)\in {\it \overline{spec}}\,
(H)$;

\item[{\bf (Shift property)}]\ $\bar{c} (a, H+\lambda (t) ) =
\bar{c} (a, H) + \int_0^1 \lambda (t) \ dt$ for any Hamiltonian
$H$ and function $\lambda:S^1 \to \R$;

\item[{\bf (Monotonicity)}]\ If $H_1\leq H_2$, then $\bar{c}
(a, H_1)\leq \bar{c} (a, H_2)$;

\item[{\bf (Lipschitz property)}]\
The map $H\mapsto \bar{c} (a, H)$ is Lipschitz on the space of
(time-dependent) Hamiltonians $H: M\times S^1\to\R$ with respect
to the $C^0$-norm;

\item[{\bf (Symplectic invariance)}]\
$\bar{c} (a,\phi^*H) = \bar{c} (a,H)$ for every $\phi \in \Symp_0
(M)$, $H \in C^{\infty} (M)$;

\item[{\bf (Normalization)}]\ $\bar{c} (a,0) = 0$ for
every even-dimensional singular homology class $a \in H_* (M,
\C)$;

\item[{\bf (Homotopy invariance)}]\ $\bar{c} (a, H_1) = \bar{c} (a, H_2)$
for any {\it normalized} $H_1, H_2$ generating the same
$\phi\in\tHam (M)$. Thus one can define $c (a,\phi)$ for any
$\phi\in\tHam (M)$ as $\bar{c} (a,H)$ for any normalized $H$
generating $\phi$. Note that $c (a, \phi)\in {\it spec}\, (\phi)$.

\item[{\bf (Triangle inequality)}]\
$c (a\ast b, \phi\psi)\leq c (a, \phi) + c (b, \psi)$.

\end{description}

\section{From a Calabi quasi-morphism to a symplectic quasi-state}
\label{sec-calabi-quasi}

In this section we prove Theorem~\ref{thm-sympl-quasistate}.
Assume that $M$ is spherically mo\-no\-tone. In this case the
Novikov ring of $M$ is a field of complex Laurent series in one
variable. The even-degree part $QH_{ev} (M)$ of $QH_* (M)$ is a
commutative algebra over this field. Assume that the algebra
$QH_{ev} (M)$ is semi-simple in the sense of \cite{EP} -- this
holds, for instance, if $M$ is one of the symplectic manifolds
listed in the statement of Theorem~\ref{thm-sympl-quasistate}: the
standard $\C P^n$, a complex Grassmannian, $\C
P^{n_1}\times\ldots\times\C P^{n_k}$ with a monotone product
symplectic structure, the monotone symplectic blow-up of $\C P^2$
at one point. Denote by ${\rm vol}\, (M^{2n}) := \int_M \omega^n $
the total symplectic volume of $M$. The main result of \cite{EP}
states that for a suitable choice of an idempotent $a \in QH_{ev}
(M)$, the function $\mu: \tHam (M)\to\R$ given by
\begin{equation}
\label{eqn-mu-def} \mu (\phi) := - {\hbox{\rm vol}\, (M)}\cdot
\lim_{k\to +\infty} c (a,\phi^k )/k
\end{equation}
is a {\it homogeneous quasi-morphism} on the group $\tHam (M)$
with a number of additional properties. More precisely, the
following holds:

\begin{itemize}

\item{\bf (Quasi-additivity)} There
exists $K>0$, which depends only on $\mu$, so that
\[
|\mu (\phi\psi) - \mu (\phi) - \mu(\psi)|\leq K \ \ {\rm for}\
{\it all}\ {\rm elements}\ \phi, \psi\in\tHam (M);
\]

\item{\bf (Homogeneity)}
$\mu (\phi^m) = m\mu (\phi)$ for each $\phi$ and each $m\in\Z$.

\end{itemize}
\medskip
\noindent

To proceed with properties of $\mu$ we need the following
notations. For a (time-dependent) Hamiltonian $H$ on $M$ write
$\phi_H$ for the element of $\tHam (M)$ represented the
identity-based path in $\Ham (M)$ given by the $[0,1]$-time
Hamiltonian flow generated by $H$. For an open $U\subset M$ denote
by $\tHam (U)\subset \tHam (M)$ the subgroup of elements generated
by Hamiltonians $H(x,t)=H_t(x)$ with ${\it supp}\, H_t\subset U$
for all $t\in S^1$. Denote by $\Cal : \tHam (U) \to\R$ the
classical Calabi homomorphism: $\Cal (\phi_H) := \int_0^1 \int_U
H_t \omega^n dt$, where ${\it supp}\, H_t\subset U$ for all $t$.
\medskip
\noindent

\medskip
\noindent
\begin{itemize}

\item{ \bf (Calabi property)} If $U\subset M$ is open
and displaceable, then the restriction of $\mu$ on $\tHam (U)
\subseteq \tHam (M)$ is the Calabi homomorphism $\Cal : \tHam (U)
\to\R$.

\item{ \bf (Lipschitz property)} $|\mu(\phi_F) -\mu(\phi_H)| \leq
{\hbox{\rm vol}\, (M)}\cdot ||F-H||_{C^{0}}.$

\end{itemize}

\medskip
\noindent Define now $\zeta:C^{\infty}(M) \to \R$ by

\begin{equation}
\label{eqn-zeta-mu} \zeta (F) = \frac{\int_M F\omega^n}{{\rm
vol}\, (M)} - \frac{\mu (\phi_F)}{{\rm vol}\, (M)} = \lim_{k \to
+\infty} \frac{\bar{c}(a,kF)}{k} \;.
\end{equation}

\medskip
\noindent Using the Lipschitz property of $\mu$, we readily extend
$\zeta$ to a functional on $C(M)$. Let us check that $\zeta$
satisfies the axioms of a symplectic quasi-state. Since $\mu(\id)
= 0$ in view of homogeneity of $\mu$, we get the normalization
axiom. Invariance and monotonicity of spectral invariants yield
the invariance and the monotonicity axioms respectively. The
Calabi property of $\mu$ yields the vanishing axiom. To check the
strong quasi-additivity axiom, note that $\phi_F$ and $\phi_G$
commute if $\{F,G\}=0$. The desired result follows from the
following general fact (which is an easy exercise): restriction of
a homogeneous quasi-morphism to any abelian subgroup is a
homomorphism. This completes the proof of
Theorem~\ref{thm-sympl-quasistate}. \Qed

\section{A partial symplectic quasi-state}
\label{sec-existence-part-qstate}

Let $M$ be a closed strongly semi-positive and rational symplectic
manifold. For an element $\phi \in \tHam (M)$ write for brevity $c
(\phi) = c ([M],\phi)$ and, as above, define $\mu$ as a
homogenization of $c ([M],\cdot)$:
\begin{equation}
\label{eqn-mu-def-1} \mu (\phi) := - {\hbox{\rm vol}\, (M)}\cdot
\lim_{k\to +\infty} c (\phi^k )/k\;.
\end{equation}
It is easy to see that $\mu$ is {\bf not} a quasi-morphism already
when $M$ is the 2-torus -- see the discussion following
Question~\ref{ques-surf-3} in Section~\ref{sec-surfaces}.
Moreover, a similar argument actually shows that for any (strongly
semi-positive, rational) symplectic direct product $M\times
\T^{2n}$ the homogenization of {\it any} spectral number $c (a,
\cdot)$ cannot be a quasi-morphism.

In spite of this, $\mu$ has a number of nice properties which will
enable us to show that the functional $\zeta$ given by
\eqref{eqn-zeta-mu} is a partial symplectic quasi-state. We shall
need the following definition. Given a displaceable open set $U
\subset M$, each $\phi\in \tHam (M)$ can be represented as a
product of elements of the form $\psi\theta\psi^{-1}$ with $\theta
\in \tHam(U)$. This follows from Banyaga's fragmentation lemma
\cite{Ban}. Denote by $\| \phi\|_U$ the minimal number of factors
in such a product.

\begin{thm}
\label{thm-partial-quasi-objects-1} Suppose $M$ is strongly
semi-positive and rational. The functional $\mu: \tHam (M)\to \R$,
given by (\ref{eqn-mu-def-1}), is well defined and has the
following properties:
\begin{itemize}
\item{\bf (Controlled quasi-additivity)} Given a
displaceable open subset $U$ of $M$, there exists a constant $K$,
depending only on $U$, so that
\[
| \mu (\phi\psi) - \mu (\phi) - \mu (\psi)| \leq K \min\, \{
\|\phi\|_U, \|\psi\|_U\}
\]
for any $\phi,\psi\in \tHam (M)$.
\item{\bf (Semi-homogeneity)}
$\mu (\phi^m) = m\mu (\phi)$ for any $\phi$ and any $m\in\Z_{\geq
0}$.
\end{itemize}
In addition it has the Calabi and Lipschitz properties defined in
the previous section.
\end{thm}

\noindent Postponing the proof of we first prove
Theorem~\ref{thm-part-quas}.

\medskip
\noindent {\bf Proof of Theorem~\ref{thm-part-quas} (assuming
Theorem~\ref{thm-partial-quasi-objects-1}):} Define a functional
$\zeta: C^{\infty}(M) \to \R$ by formula \eqref{eqn-zeta-mu}. We
claim that $\zeta$ is a partial symplectic quasi-state. Arguing
exactly as in the end of the previous section, we check the
monotonicity, vanishing, normalization and invariance axioms.
Semi-homogeneity of $\mu$ yields that $\zeta(\lambda F) = \lambda
\zeta(F)$ for $\lambda \in \N$ and all smooth $F$. As a logical
consequence we get that the same holds for all positive rational
$\lambda$. Using the Lipschitz property of $\mu$, we pass to the
limit and get this for all positive $\lambda$, thus establishing
the semi-homogeneity axiom.

It remains to verify the partial additivity axiom. Assume that
$\{F,H\}= 0$ and $\supp\, H$ is contained in a displaceable open
subset $U$. Note that $||\phi_H^k||_{U}=1$ for all $k \in \N$.
Since $\phi_F$ and $\phi_H$ commute we have (using controlled
quasi-additivity of $\mu$)
$$\mu(\phi_F\phi_H) = \frac{1}{k} \mu((\phi_F\phi_H)^k) =
\frac{1}{k}(k\mu(\phi_F) + k\mu(\phi_H) + r_k)\;,$$ where
$|r_k|\leq K$. Taking the limit as $k \to +\infty$ we get that
$$\mu(\phi_F\phi_H) = \mu(\phi_F)+\mu(\phi_H) = \mu(\phi_F)\ + \int_M H
\omega^n\;,$$ where the last equality follows from the Calabi
property and the fact that $\supp\, H$ is displaceable. Further,
$\phi_{F+H} = \phi_F\phi_H $ since $F$ and $H$ commute.
Substituting this into the definition of $\zeta$, we get
$\zeta(F+H) = \zeta(F)$, as required. This completes the proof.
\Qed

\medskip
\noindent {\bf Proof of
Theorem~\ref{thm-partial-quasi-objects-1}:} The proof is divided
into a sequence of lemmas. In what follows we fix an open
displaceable subset $U$ of $M$ and write for simplicity
$||\phi||:= ||\phi||_U$.

\begin{lemma}
There exists a constant $C \geq 0$ such that for any $\phi\in
\tHam (U) $
\[
0\leq c ( \phi) + c ( \phi^{-1})\leq C .
\]
\end{lemma}

\smallskip
\noindent {\bf Proof.} Suppose $f\in \tHam (M)$ is a lift of a
Hamiltonian diffeomorphism displacing $U$. Then the "shift of the
spectrum" trick of Y.Ostrover \cite{Ostr1} (cf. \cite{En-comm},
\cite{EP}) yields that for a certain $E\in\R$, depending on
$\phi$,
\[
c ( f \phi) = c ( f) + E,
\]
\[
c ( f \phi^{-1}) = c ( f) - E.
\]
Here we use the spectrality and the Lipschitz property of $c$. The
signs in the formulae above comply with the sign convention as in
\cite{EP}. Thus
\[
c ( f \phi) + c ( f \phi^{-1}) = 2 c (f).
\]
In view of the triangle inequality\footnote{Note that, since $[M]$
is the unit in $QH_\ast (M)$, the triangle inequality for
$c(\cdot) = c ([M], \cdot )$ has the form $c (\phi\psi)\leq
c(\phi) + c(\psi)$.},
\[
0\leq c ( \phi) + c ( \phi^{-1}),
\]
\[
c ( \phi) \leq c ( f \phi ) + c ( f^{-1}),
\]
\[
c ( \phi^{-1}) \leq c ( f \phi^{-1} ) + c ( f^{-1}).
\]
Hence
\[
0\leq c ( \phi) + c ( \phi^{-1}) \leq c ( f \phi) + c ( f
\phi^{-1}) + 2 c ( f^{-1})\leq 2 c ( f) + 2 c (f^{-1}).
\]
Set $C := 2c ( f) + 2c ( f^{-1})$. This is a non-negative number
because of the triangle inequality. The lemma is proved. \Qed

\begin{lemma}
\label{lemma-q} For any $\phi\in \tHam (U)$ and any $\psi\in \tHam
(M)$ one has
\[
c ( \phi) + c ( \psi) - C\leq c ( \phi \psi)\leq c ( \phi) + c (
\psi),
\]
where $C$ is the constant from the previous lemma.
\end{lemma}

\smallskip
\noindent {\bf Proof.} The second inequality is just the triangle
inequality. To obtain the first one, observe that the triangle
inequality yields
\[
c ( \psi) \leq c ( \phi \psi) + c ( \phi^{-1}).
\]
This, along with the previous lemma, implies
\[
c ( \phi \psi) \geq c ( \psi) - c ( \phi^{-1})\geq c ( \psi) + c (
\phi) - C.
\]
\Qed

Using a straightforward inductive argument one generalizes the
lemma above as follows. Take any $\phi_1,\ldots , \phi_m, \psi \in
\tHam (M)$ with $||\phi_i||=1$ for all $i$. Then
\begin{equation}
\label{eqn-aa1} | c (\phi_1\cdot\ldots\cdot \phi_m\psi) -
\sum_{i=1}^m c (\phi_i) -c(\psi) | \leq mC.
\end{equation}

This formula (with $\psi=\id$) yields
\begin{equation}
\label{eqn-aa2} | c ((\phi_1\cdot\ldots\cdot \phi_m)^l) -
l\sum_{i=1}^m c (\phi_i) | \leq l mC.
\end{equation}

Take any $\phi \in \Ham(M)$ and represent it as
$\phi=\phi_1\cdot\ldots\cdot\phi_m$ with $||\phi_i||=1$ for all
$i$. Formula \eqref{eqn-aa2} implies that for some large enough
positive $E$ (depending on $\phi$) the sequence $ \{ c (\phi^l) +E
l\}_{l\in\N}$ is non-negative. On the other hand, because of the
triangle inequality, this sequence is sub-additive. This yields
the existence and finiteness of $\lim_{l\to +\infty} (c (\phi^l)
+E l)/l$ and, accordingly, of $\lim_{l\to +\infty} c (\phi^l)/l$.
Therefore the function $\mu$ is well defined. The semi-homogeneity
of $\mu$ follows immediately from its definition. The proof of the
Lipschitz property of $\mu$ simply repeats the proof of a similar
Proposition 3.5 in \cite{EP}.

\medskip
\noindent Now we are going to check controlled quasi-additivity of
$\mu$. Assume without loss of generality that the volume of $M$
equals $1$, so that $$\mu(\phi) = -\lim_{k \to +\infty}
c(\phi^k)/k\;.$$ We claim that for $\phi,\psi \neq \id$,
\begin{equation}\label{eq-contr-1}
| \mu (\phi\psi)-\mu (\phi)-\mu (\psi) | \leq 2C\cdot \min
(2||\phi||-1, 2 ||\psi||-1).
\end{equation}
The controlled quasi-additivity follows immediately from
(\ref{eq-contr-1}) if one sets $K:= 4 C$. We prove the claim by
induction on $m:= \min(||\phi||,||\psi||)$.

\medskip
\noindent{\sc Induction basis $m=1$:} Assume without loss of
generality that $||\phi||=1$. Note that $$(\phi\psi)^k = \bigg(
\prod_{i=0}^{i=k-1} \psi^i\phi\psi^{-i}\bigg)\cdot \psi^k.$$
Applying \eqref{eqn-aa1} and using the conjugation invariance of
$c (\cdot)$ we get $$|c((\phi\psi)^k) -kc(\phi) -c(\psi^k)| \leq
Ck.$$ Combining this with inequality $$|c(\phi^k) -kc(\phi)| \leq
Ck,$$ which follows from \eqref{eqn-aa2}, dividing by $k$ and
passing to the limit as $k \to +\infty$ we get the desired result.

\medskip
\noindent {\sc Induction step $m \mapsto m+1$:} Assume without
loss of generality that $||\phi|| = m+1$. Then $\phi$ can be
decomposed as $\phi= \phi_m\phi_1$ where $||\phi_m||=m$ and
$||\phi_1||=1$. Using the induction assumption we have
$$|\mu(\phi_m\phi_1\psi) -\mu(\phi_m) -\mu(\phi_1\psi)| \leq
2C(2m-1),$$ $$|\mu(\phi_1\psi) -\mu(\phi_1) -\mu(\psi)| \leq 2C $$
and $$|\mu(\phi_1)+\mu(\phi_m) -\mu(\phi_m\phi_1)| \leq 2C.$$
Adding up these inequalities we get that $$|\mu(\phi\psi)
-\mu(\phi) -\mu(\psi)| \leq 2C(2m+1),$$ as desired. This completes
the proof of the claim and of the controlled quasi-addi\-ti\-vity.

Finally, the proof of the Calabi property of $\mu$ virtually
repeats the proof of a similar Proposition 3.3 in \cite{EP}. The
symplectic invariance of $\mu$ follows from the symplectic
invariance of the spectral numbers. This finishes the proof of
Theorem~\ref{thm-partial-quasi-objects-1}. \Qed

\section{Symplectic quasi-states on surfaces}
\label{sec-surfaces}

\medskip
\noindent {\sc Symplectic quasi-measures.} A quasi-measure on a
symplectic manifold $M$ is called {\it symplectic} if it is
$\Symp_0 (M)$-invariant and vanishes on displaceable closed
subsets.

Here we discuss this notion in the case when $M$ is a closed
surface equip\-ped with an area form. According to the general
construction from \cite{Aar91}, any quasi-measure $\tau$ gives
rise to a quasi-state $\zeta_\tau$. Roughly speaking, the
definition of $\zeta_\tau$ is as follows. For a function $F \in
C(M)$ define a measure $\sigma_F$ on $\R$ by its values on
intervals
$$\sigma_F ([a;b)):= \tau(\{F \geq a\})-\tau(\{F \geq b\},$$
and put $\zeta_{\tau}(F) := \int_{\R} s\cdot d\sigma_F(s)$. If
$\tau$ is a symplectic quasi-measure, the quasi-state
$\zeta_{\tau}$ automatically satisfies all the axioms of a
symplectic quasi-state except, probably, strong quasi-additivity
stating that $\zeta_{\tau}$ is linear on the centralizer (with
respect to the Poisson bracket) of any smooth function $F$.

\begin{thm}
\label{surf} On a closed surface, the strong quasi-additivity
axiom follows from the usual quasi-linearity. In particular, any
symplectic quasi-measure gives rise to a symplectic quasi-state.
\end{thm}

\medskip
\noindent {\bf Proof:} Let $F,G$ be a pair of $C^{\infty}$-smooth
functions on a closed surface $M$ with $\{F,G\} = 0$. The
Poisson-commutativity can be interpreted as follows: the
differential of the map $$\Phi: M \to \R^2,\; x \mapsto
(F(x),G(x)),$$ has rank $\leq 1$ for at each point $x \in M$. Put
$\Delta:= \text{Image} (\Phi)$. Denote by $d_c$ and $d_{h}$ the
covering dimension and the Hausdorff dimension of $\Delta$
respectively. It is a standard fact of dimension theory that $d_c
\leq d_{h}$, see e.g. the proof of Theorem (6.2.10) in Edgar's
book \cite{Edgar}. Further, $d_h \leq 1$. This follows from a
result of Dubovickii \cite{Dub} which is a partial case of a more
general theorem of Sard \cite{Sard}. Therefore $d_c \leq 1$.

Define a quasi-state $\eta$ on $C(\Delta)$ by $\eta (H) :=
\zeta_{\tau}(\Phi^*H)$. The Wheeler-Shakh\-ma\-tov Theorem
\cite{Grubb-LaBerge}, \cite{Wheeler} implies that every
quasi-state on a normal topological space (and hence on any metric
space) of covering dimension $\leq 1$ is linear. Hence $\eta$ is
linear. Applying this result to the restriction of the coordinate
functions on $\R^2$ to $\Delta$ we get that
\begin{equation}
\label{eq-surf-add} \zeta_{\tau}(F+G) = \zeta_{\tau}(F)
+\zeta_{\tau}(G), \end{equation} as required. \Qed

\medskip
\noindent

Note that in the proof above we used that the functions $F$ and
$G$ are infinitely smooth in order to deduce inequality $d_h \leq
1$ from the Dubovickii-Sard theorem.

\begin{problem}\label{prob-surf-1}
Extend identity \eqref{eq-surf-add} to Poisson-commuting functions
of finite smoothness.
\end{problem}

\medskip
\noindent For instance, one can try to find a uniform
approximation of the pair $(F,G)$ by a Poisson-commuting pair of
$C^{\infty}$-functions.

\medskip
\begin{rem} \label{quest-surf-1}{\rm
In contrast to the case of surfaces,  the only known to us example
of a symplectic quasi-measure on higher-dimensional manifolds
comes from the ``Floer-homological" symplectic quasi-state whose
existence is established in Theorem \ref{thm-sympl-quasistate}.}
\end{rem}

\medskip
\begin{rem}
\label{rem-non-sympl-qstates-higher-dim} {\rm In the case ${\rm
dim}\ M > 2$ D.Grubb \cite{Grubb-private} constructed examples of
quasi-states which are not strongly quasi-additive. The
quasi-measures in the examples of Grubb do not necessarily vanish
on displaceable sets (and hence are not symplectic) but may vanish
on a stem. The push-forward of such a quasi-measure by a moment
map of a finite-dimensional Poisson-commutative subspace of
$C^\infty (M)$ is not necessarily a measure. }
\end{rem}

\medskip

Now we address a question about existence and uniqueness of
symplectic quasi-states and quasi-measures on surfaces.

\medskip
\noindent {\sc The 2-sphere.} The group $\Ham (S^2)$ admits a
Calabi quasi-morphism \cite{EP}, which in accordance with our
discussion in Section~\ref{sec-calabi-quasi} yields existence of a
symplectic quasi-state and a symplectic quasi-measure on $C
(S^2)$. Theorem 5.2 in \cite{EP} shows that any two Calabi
quasi-morphisms on $\Ham (S^2)$ coincide on the set of elements
generated by time-independent Hamiltonians. The same argument
proves that any two symplectic quasi-states coincide on the set of
smooth Morse functions on $S^2$. Hence $C (S^2)$ {\it carries
unique symplectic quasi-state and quasi-measure}.

An explicit calculation presented in \cite{EP} shows that the
restriction of this symplectic quasi-state, say $\zeta$, to the
subalgebra $\cA_F \subset C(M)$ generated by a single Morse
function $F \in C^{\infty} (S^2)$ is multiplicative: $\zeta (GH) =
\zeta(G)\zeta(H)$ for all $G,H \in \cA_F$. Using this along with
the continuity of $\zeta$ one can easily show that $\zeta$ is
multiplicative on $\cA_F$ for {\it any} $F\in C (S^2)$. Now a
theorem of Aarnes \cite{Aar-MathAnn} yields that {\it the
corresponding quasi-measure is simple: it takes values $0$ and $1$
only.} It is unclear whether this phenomenon persists in higher
dimensions, thus we pose the next question.

\begin{question}\label{ques-surf-10}
Consider the ``Floer-homological" symplectic quasi-state $\zeta$
on the complex projective space $\C P^n$ constructed in Theorem
\ref{thm-sympl-quasistate}. Is it multiplicative when $n \geq 2$?
In particular, is it true that $\zeta(F^2)=\zeta(F)^2$ for all
continuous functions $F$ on $\C P^n$ ?
\end{question}

For completeness, we present the formula for $\zeta$ on $\cA_F$,
where $F$ is a Morse function, obtained in \cite{EP}. Assume that
the total area of the sphere equals 1. One shows that there exists
unique (may be, singular) connected component of a level set of
$F$, say $\gamma$, so that the area of any connected component of
$S^2 \setminus \gamma$ is $\leq \frac{1}{2}$. Note that every $G
\in \cA_F$ is constant on connected components of level sets of
$F$. It turns out that
$$\zeta(G) = G(\gamma).$$

A symplectic quasi-measure $\tau$ corresponding to $\zeta$ can be
described as follows (we thank D.Grubb who pointed this to us). A
set $A\subset S^2$ is called {\it solid} if both $A$ and
$S^2\setminus A$ are connected. According to the results of Aarnes
\cite{Aar-solid-sets} and Aarnes and Rustad \cite{Aar-Rus}, the
quasi-measure $\tau$ is completely defined by the following
condition: for a closed solid set $A\subset S^2$ one has $\tau (A)
= 1$ if the Lebesgue measure of $A$ is greater or equal to $1/2$
and $\tau (A) = 0$ otherwise.

\medskip
\noindent {\sc The 2-torus.} Existence of a symplectic
quasi-measure, say $\tau$, in this case follows from a work by
Grubb (see Theorem 32 of \cite{Grubb-TAMS}, where the auxiliary
quasi-measures used in the definition of $\tau$ are taken to be
the standard Lebesgue measure). The value of $\tau$ on any
$2$-dimensional smooth connected closed submanifold with boundary
$W \subset \T^2$ can be calculated as follows (see Theorem 32 of
\cite{Grubb-TAMS}). If $W$ is contractible in $\T^2$ we have
$\tau(W)=0$. If $W$ is non-contractible and $\partial W$ has
$k\geq 0$ contractible connected components that bound pair-wise
disjoint discs $D_1,...,D_k$ (in case $k=0$ there are no discs),
then
$$\tau(W) = \text{Area}(W) + \sum_{i=1}^k \text{Area}(D_i).$$

\medskip
\begin{rem}\label{ques-surf-2}{\rm
It would be interesting to describe all symplectic
quasi-measu\-res on the 2-torus; for more examples of such
quasi-measures see a recent preprint \cite{Knudsen} by Knudsen. }
\end{rem}

\medskip
\noindent By Theorem~\ref{surf} above, a symplectic quasi-measure
on $T^2$ gives rise to a symplectic quasi-state.

\begin{question} \label{ques-surf-3}
Is Grubb's symplectic quasi-measure associated to a quasi-morphism
on $\Ham(\T^2)$?
\end{question}

Such a quasi-morphism, if exists, cannot come from spectral
numbers described in Section~\ref{sec-spectral-numbers}. To see
this denote by $\tau$ any symplectic quasi-measure on $\T^2$.
Introduce coordinates $(p,q) \mod \;1$ on $\T^2$ so that the
symplectic form is given by $dp\wedge dq$. Let $\alpha = \{p=0\}$
and $\beta = \{p=1/2\}$ be two meridians dividing the torus into
two open annuli $A = \{p \in (0;1/2)\}$ and $B= \{p \in (1/2;1)\}$
of equal area. Note that
$$\tau(A)+\tau(B)+\tau(\alpha)+\tau(\beta) =1.$$ The
$\Symp_0$-invariance of $\tau$ yields $\tau(A)=\tau(B)$ as well as
$\tau(\alpha) = \tau(\beta) = 0$ (the torus contains an
arbitrarily large number of pair-wise disjoint symplectic shifts
of a meridian). Thus, putting $A' = A \cup \alpha \cup \beta$, we
have $\tau(A') = 1/2$. On the other hand, choose a sequence of
cut-off functions $F_{i}(p)$ approximating the characteristic
function of $A'$ so that the only critical values of $F_i$ are $0$
and $1$. The key feature of the Hamiltonian flow generated by
$F_i$ is that its only {\it contractible} closed orbits are the
critical points, hence the action spectrum ${\it
\overline{spec}}\, (tF_i)$ equals $\{0;t\}$. Hence, using
continuous dependence of spectral numbers on the Hamiltonian, we
get that for every homology class $a \in H_*(\T^2)$, we have
either $\bar{c} (a,tF_i) = 0$ or $\bar{c} (a,tF_i) = t$.
Substituting this into the right term of formula
\eqref{eqn-zeta-mu}, we get that $\zeta(F_i)$, if well defined,
must be either $0$ or $1$ and hence $\tau(A') \neq 1/2$. This
contradiction proves the claim.

Note that Hamiltonians $F_i$ above have a wealth of
non-contractible periodic orbits. In principle, the symplectic
field theory \cite{EGH}, or, more precisely, its version called
branched Floer homology (work in progress by V. Ginzburg and E.
Kerman) which deals with Hamiltonian diffeomorphisms, may lead to
a generalization of spectral numbers which takes into account
non-contractible orbits as well. It would be interesting to
understand whether this path leads to a symplectic quasi-measure.

\section{Digging out a stem} \label{sec-stem}

Assume that one faces the problem of the following type: "Prove
that a certain specific Lagrangian submanifold $L$ of a symplectic
manifold $M$ is non-displac\-eable". The mainstream approach to
this problem is to show that the Lagrangian Floer homology of $L$
is well defined and does not vanish. Our results above give rise
to another potential approach (cf. \cite{BEP}): show that $L$ is a
stem (see Definition~\ref{def-stem}) and deduce the
non-displaceability of $L$ from
Theorem~\ref{thm-commut-alg-has-nondispl-fiber}. Let us emphasize
that this approach is not "soft": unveiling the proof, one sees
that we use an information about the asymptotic behaviour of
Hamiltonian Floer homology for Hamiltonians concentrated near $L$.
While in certain situations our method is simpler, it does not
provide a lower bound on the number of intersections (assuming
they are all transversal) between $L$ and its image under a
Hamiltonian isotopy -- a bound which is usually given by the
Lagrangian Floer homology approach whenever it works.

Below we illustrate our approach for the Lagrangian Clifford torus
in $\C P^n$ and for a similar torus in a monotone blow-up of $\C
P^2$ at one point.

\bigskip
\noindent {\sc The Clifford torus in $\C P^n$.} This example is
taken from \cite{BEP}. Let $M$ be $\C P^n$ with the Fubini-Study
symplectic form. Consider the standard Hamiltonian $\T^n$-action
on $M$ whose moment polytope is a simplex in $\R^n$. Denote by $L$
the Lagrangian torus which is the fiber of the moment map over the
barycenter of the simplex -- it is called the {\it Clifford torus}
and can be described as
\[
L := \{\, [z_0:\ldots : z_n]\in \C P^n \ | \ \ |z_0| =\ldots =
|z_n|\, \}.
\]
All the fibers of the moment map, other than $L$, are displaceable
-- this easily follows from the observation that permutations of
homogeneous coordinates can be realized by Hamiltonian
diffeomorphisms of $\C P^n$ coming from the natural action of $PU
(n+1)$ on $\C P^n$. Thus $L$ is a stem. Hence, according to
Corollary~\ref{cor-toric}, $L$ is non-displaceable \cite{BEP}.

In fact, the non-displaceability of $L$ can be also proved by
means of the Lagrangian Floer homology. In this way Cho \cite{Cho}
showed that if an image of $L$ under a Hamiltonian isotopy is
transversal to $L$ then the number of intersections between them
must be at least $2^n$ (which is the sum of the Betti numbers of
$L$).

\bigskip \noindent {\sc The Clifford torus in the monotone blow-up of
$\C P^2$ at one point.} Our interest in this example is due to the
fact that in this case the obstructions to displaceability of $L$
coming from the Lagrangian Floer homology do vanish according to a
result by Cho and Oh \cite{Cho-Oh}.

Here is the description of $L$. Consider a spherical shell $W$
lying in the standard symplectic linear space $\C^2$: $$W = \{\,
(u_1,u_2) \in \C^2\ \big| \ \frac{1}{3} \leq \pi(|u_1|^2+|u_2|^2)
\leq 1\, \}.$$ Making a symplectic cut (i.e. collapsing the
boundaries along the fibers of the characteristic foliation) we
get a closed symplectic 4-manifold $M$ which is one of the models
of the blow up of $\C P^2$ at one point. The details of the
construction can be extracted from the description of the
symplectic structure on $M$ given in \cite{Aud}, page 61. The
symplectic manifold $M$ is spherically monotone. A Lagrangian
torus $$L := \{\pi |u_1|^2 = \pi |u_2|^2 = \frac{1}{3}\}$$ is
called {\it the Clifford torus} of $M$ -- it can be viewed as the
Clifford torus in $\C P^2$ which "survived" the blow-up.

\begin{thm}\label{thm-clif}
The Clifford torus $L \subset M$ is a stem.
\end{thm}

\medskip
\noindent Combining this with Corollary~\ref{cor-toric} we get
that $L$ is non-displaceable.

\medskip
\noindent {\bf Proof:} Consider a Hamiltonian action of the
2-torus on $M$, which in the spherical shell model is defined by
its moment map
$$\Phi: W \to \R^2\;,\; (u_1,u_2) \mapsto (\pi |u_1|^2,\pi
|u_2|^2)\;.$$ We shall show that $L$ is the stem of a
Poisson-commutative subspace generated by the coordinate functions
of $\Phi$. The image $\Delta$ of $\Phi$ is a trapezoid $ABCD$ in
the plane with the vertices $$A = (0,1/3), B = (1/3,0), C =
(1,0),D = (0,1).$$ The Clifford torus $L$ is given by
$\Phi^{-1}(Q),$ where $Q = (1/3,1/3)$.

\medskip
\noindent {\sc Claim:} The fiber $\Phi^{-1}(X)$ is displaceable
for every $X \neq Q$.

\medskip
\noindent We use the following notation for lines and segments on
the plane: $PR$ stands for the line passing through points $P$ and
$R$,$\;$ $[PR)$ denotes the segment with vertices $P$ and $R$ so
that $P$ is included and $R$ is excluded and so on. We write
$|PR|$ for the Euclidean length of $[PR]$.

\medskip Consider the points $$P = (1/6,1/6) \in [AB],\;R=(1/2,1/2) \in
[CD]\;.$$

\medskip
\noindent {\sc Case I: $X \notin [PR]$.} The unitary
transformation $S:(u_1,u_2) \to (u_2,u_1)$ of $W$ commutes with
the $\T^2$-action and induces the symmetry of $\Delta$ over the
line $PR$ which sends $X$ to a point $X' \neq X$. Hence
$S(\Phi^{-1}(X)) \cap \Phi^{-1}(X) = \emptyset$, which proves the
claim in this case.

\medskip
\noindent In order to proceed further, take the point $E=
(2/3,1/3) \in [CD]$. The segment $[AE]$ divides $\Delta$ into a
triangle $\Delta'$ and a parallelogram $\Pi$. We assume that
$\Delta'$ contains segments $[DA)$ and $[DE)$ and does not contain
$[AE]$, while $\Pi$ contains $[BA)$ and $[BC)$ and does not
contain the two other edges.

\medskip
\noindent {\sc Case II: $X \in (QR]$.} The set
$\Phi^{-1}(\Delta')$ is $\T^2$-equivariantly symplectomorphic to
the standard symplectic ball $\{\pi(|w_1|^2+|w_2|^2) < 2/3\}$ (the
vertex $D$ corresponds to the center of the ball). This follows
from the local version of Delzant theorem -- see \cite{Karshon}.
The unitary transformation $(w_1,w_2) \to (w_2,w_1)$ of this ball
commutes with the $\T^2$-action and induces an affine involution
of $\Delta'$ whose fixed point set coincides with $[DQ)$. This
involution sends $X$ to some point $X' \neq X$. We conclude that
the torus $\Phi^{-1}(X)$ can be sent to $\Phi^{-1}(X')$ by a
Hamiltonian isotopy, and therefore is displaceable.

\medskip
\noindent {\sc Case III: $X \in [PQ)$.} The set $\Phi^{-1}(\Pi)$
is $\T^2$-equivariantly symplectomorphic to the standard
symplectic polydisc $\cC = \cD _1 \times \cD _2$ with $$\cD_1=
\{\pi|w_1|^2 < 2/3\},\; \cD_2= \{\pi|w_2|^2 < 1/3\}$$ (the vertex
$B$ corresponds to the center of the polydisc). This again follows
from the local version of Delzant theorem \cite{Karshon}. The
projection of $\cC$ to $\cD _1$ sends the torus $\Phi^{-1}(X)$ to
a circle $\Gamma:=\{\pi |w_1|^2 = r\}$ which encloses a disc of
area $r$. The area $r$ corresponding to the point $X$ can be
calculated as follows. Let $Y$ be the projection of $X$ to $BC$
along $AB$. Then
$$\frac{r}{\text{Area}(\cD _1)} = \frac{|BY|}{|BC|} <
\frac{1}{2}\;.$$ This inequality guarantees that $\Gamma$ is
displaceable in $\cD _1$ by a Hamiltonian transformation of $\cD
_1$. Lifting this transformation to $\cC$ we get that
$\Phi^{-1}(X)$ is displaceable. This completes the proof of the
claim, and hence of the theorem. \Qed

\medskip
\noindent Sometimes even in seemingly simple situations it is hard
to decide whether a given Lagrangian submanifold is a stem. For
instance, we do not know an answer to the following question:

\begin{question} \label{ques-stem}
Consider $\R P^2 \subset \C P^2$ or the anti-diagonal in the
monotone $S^2 \times S^2$. Are they stems?
\end{question}

\section{On the history and the physical meaning of quasi-states}
\label{sec-hist}

The notion of quasi-state has an amusing history. To discuss it
let us recall the mathematical model of quantum mechanics which
goes back to von Neumann's famous book \cite{von Neumann}
published in 1932 : Its basic ingredients are the real Lie algebra
of observables $\cA_q$ ($q$ for quantum) whose elements (in the
simplest version of the theory) are hermitian operators on a
finite-dimensional complex Hilbert space $H$ and the Lie bracket
is given by $$[A,B]_{\hbar} = \frac{i}{\hbar}(AB-BA)\;,$$ where
$\hbar$ is the Planck constant. Observables represent physical
quantities such as energy, position, momentum etc. The state of a
quantum system is given by a functional $\zeta: \cA_q \to \R$
which satisfies the following axioms:

\medskip
\noindent {\bf (Additivity)} $\zeta(A+B) = \zeta(A)+\zeta(B)$ for
all $A,B \in \cA_q$;

\medskip
\noindent {\bf (Homogeneity)} $\zeta(cA) = c\zeta(A)$ for all $c
\in \R$ and $A \in \cA_q$;

\medskip
\noindent{\bf (Positivity)} $\zeta(A)\geq 0$ provided $A \geq 0$;

\medskip
\noindent {\bf (Normalization)} $\zeta(\id) = 1$.

\medskip

As a consequence of these axioms von Neumann proved that for every
quantum state $\zeta$ there exists a non-negative Hermitian
operator $U_\zeta$ with trace $1$ such that $\zeta (A) =
\text{tr}(U_{\zeta}A)$ for all $A \in \cA_q$. An easy consequence
of this formula is that for every state $\zeta$ there exists an
observable $A$ such that
\begin{equation} \label{eq-disp}  \zeta(A^2) -\zeta(A)^2 > 0\;.
\end{equation}
In his book von Neumann adopted a statistical interpretation of
quantum mechanics according to which the value $\zeta(A)$ is
considered as the expectation of a physical quantity represented
by $A$ in the state $\zeta$. In this interpretation the equation
\eqref{eq-disp} says that there are no dispersion-free states.
This result led von Neumann to a conclusion which in the language
of quantum mechanics can be formulated as the impossibility to
introduce hidden variables into the quantum theory. This
conclusion caused a (seemingly never ending) discussion among
physicists which (citing Ballentine \cite{Ballentine}, p. 374)
``was unfortunately clouded by emotionalism".  A number of
prominent physicists, including Bohm and Bell, disagreed with the
additivity axiom of a quantum state. Their reasoning was that the
formula $\zeta(A+B) = \zeta(A)+\zeta(B)$ makes sense {\it a
priori} only if observables $A$ and $B$ are simultaneously
measurable, that is commute: $[A,B]_{\hbar} = 0$. We refer to
Bell's paper \cite{Bell} for an account of this discussion.

In 1957 Gleason \cite{Gleason} proved a remarkable rigidity-type
theorem which can be considered as an additional argument in favor
of von Neumann's additivity axiom. Recall that two hermitian
operators on a finite-dimensional Hilbert space commute if and
only if they can be written as polynomials of the same
self-adjoint operator. Let us introduce {\it a quasi-state} on
$\cA_q$ as a real-valued functional which satisfies the
homogeneity, positivity and normalization axioms above, while the
additivity axiom is replaced by one of the two {\bf equivalent}
axioms:

\medskip
\noindent {\bf (Quasi-additivity-I)} $\zeta (A+B) =
\zeta(A)+\zeta(B)$ provided $A$ and $B$ commute: $[A,B]_{\hbar} =
0$;

\medskip
\noindent {\bf (Quasi-additivity-II)} $\zeta (A+B) =
\zeta(A)+\zeta(B)$ provided $A$ and $B$ belong to a
single-generated subalgebra of $\cA_q$.

According to the Gleason theorem, every quasi-state on $\cA_q$ is
linear (that is, a state) provided the complex dimension of the
Hilbert space $H$ is at least $3$ (it is an easy exercise to show
that in the two-dimensional case there are plenty of non-linear
quasi-states).

Let us turn now to the mathematical model of classical mechanics.
Here the algebra $\cA_c$ of observables ($c$ for classical) is the
space of continuous functions $C(M)$ on a symplectic manifold $M$.
The Lie bracket is defined as the Poisson bracket on the dense
subspace $C^{\infty}(M) \subset C(M)$. A natural question is
whether the conclusion of the Gleason theorem remains valid in the
classical context. We immediately face a dilemma: which of two
definitions of quasi-additivity one should adopt as the starting
point of such an extension. Adopting the second one, we arrive to
the definition of a quasi-state given by Aarnes. It does not
involve the symplectic structure and gives rise to the theory of
quasi-states on general topological spaces. Adopting the first
one, and taking into account the Correspondence Principle
according to which the bracket $[\;,\;]_{\hbar}$ corresponds to
the Poisson bracket $\{\;,\,\}$ in the classical limit $\hbar \to
0$, we get a definition which involves the strong quasi-additivity
axiom: $\zeta(F+G) =\zeta(F)+\zeta(G)$ whenever $\{F,G\}=0$, see
Section~\ref{sec-quasi-objects}. According to Theorem~\ref{surf}
above both definitions coincide in dimension 2. However, as it was
mentioned in Remark~\ref{rem-non-sympl-qstates-higher-dim}, strong
quasi-additivity is strictly stronger in higher dimensions. For
the sake of brevity, we refer to non-linear strongly
quasi-additive quasi-states on symplectic manifolds as to {\it
strong quasi-states}.

In light of this discussion, Theorems~\ref{thm-sympl-quasistate}
and \ref{surf} above which establish the existence of strong
quasi-states on certain symplectic manifolds can be viewed as an
"anti-Gleason phenomenon" in classical mechanics. This
interpretation is far from being transparent. Let us indicate two
points which require further clarification.

First, recall that the algebra $\cA_c$ of classical observables
can be considered as a suitable limit of matrix algebras $\cA_q$
where the dimension $N$ of the underlying Hilbert space $H$ tends
to $\infty$ and the Planck constant $\hbar$ tends to $0$. We refer
the reader to Madore's paper \cite{Madore} dealing with the case
where the classical phase space is the 2-dimensional sphere. For
certain symplectic manifolds the algebra $\cA_c$ carries a strong
quasi-state, say, $\zeta$. At the same time the Gleason theorem
rules out existence of a non-linear quasi-state on $\cA_q$ for
every given values of $N$ and $\hbar$. It would be interesting to
understand what is a footprint of $\zeta$ in the quantum world.
For instance, do the algebras $\cA_q$ carry a weaker object (a
kind of ``approximate quasi-state" still to be defined) which
converges to $\zeta$?

Second, by the analogy with quantum mechanics, one can speculate
that Poisson non-commuting functions $F$ and $G$ with  $\zeta(F+G)
\neq \zeta(F)+\zeta(G)$ are not simultaneously measurable. Does
there exist an explanation of this phenomenon in terms of
classical mechanics?  An extra difficulty here is due to the fact
that some strong quasi-states are dispersion-free. Therefore one
cannot refer to uncertainty as to the reason for the lack of
simultaneous measurability.

\bigskip
\noindent {\bf Acknowledgements.} We are grateful to P. Biran for
his generous help with this paper. We also thank D.Grubb,
Y.Karshon, Y.Ostrover, M.Syming\-ton and Y.Yomdin for useful
discussions, and the anonymous referee for corrections.

\bigskip

\bibliographystyle{alpha}

\bigskip

\noindent

\begin{tabular}{@{} l @{\ \ \ \ \ \ \ \ \ \ \,} l }
Michael Entov & Leonid Polterovich \\
Department of Mathematics & School of Mathematical Sciences \\
Technion & Tel Aviv University \\
Haifa 32000, Israel & Tel Aviv 69978, Israel \\
entov@math.technion.ac.il &
polterov@post.tau.ac.il\\
\end{tabular}

\end{document}